\documentclass[a4paper,10pt]{article}
\usepackage{paper-en}\usepackage{hyperref}

\def\thetitle{Braid spaces have no geodesic loops}
\def\theauthors{Dmitri Panov and Anton Petrunin}

\hypersetup{colorlinks=true,
citecolor=black,
linkcolor=black,
anchorcolor=black,
filecolor=black,
menucolor=black,
urlcolor=black,
pdftitle={\thetitle},
pdfauthor={\theauthors}
}

\begin{document}

\title{\thetitle}
\author{\theauthors}
\date{}
\maketitle

\begin{abstract}
The braid spaces are defined as universal metric coverings of the standard complex spaces, branched over certain hyperplane arrangements.
We describe a natural family of convex functions on braid spaces and use them to prove that braid spaces contain no geodesic loops.
This statement is a consequence of the well-known conjecture that every braid space is CAT(0).
\end{abstract}

This note expands the comments on the second problem in \cite{petrunin2026}.
It is meant to be read as a letter (it has a postscript); we wrote it with Daniel Allcock, Ruth Charney, and Michael Davis in mind, but it is addressed to any friend of $\CAT(0)$ spaces.

Let $\CC^n$ be the complex space
with coordinates $z_1,\ldots,z_n$ and the standard Euclidean metric.
Remove from $\CC^n$ the complex hyperplanes $z_i=z_j$ for all $i\ne j$;
denote the remaining set by $\WW^n$.
Let $\tilde \WW^n$ be the universal cover of $\WW^n$, equipped with the induced metric.
The completion of $\tilde \WW^n$ will be called the \emph{braid space}
with $n$ strings and denoted by $\CB^n$.

The set $\WW^n$ can be interpreted as the configuration space of $n$ punctures in the plane.
Therefore a path in $\WW^n$ describes a braid with $n$ strings, which explains our terminology.

The following conjecture is a special case of Daniel Allcock's conjecture \cite[Conjecture 1.4]{allcock}, which is related to an earlier conjecture of Ruth Charney and Michael Davis \cite[Conjecture 3]{charney-davis-95}, which in turn is motivated by the $K(\pi,1)$-conjecture \cite{vanderlek,paris}.
The authors are responsible for a more general version of Allcock's conjecture \cite[Conjecture 1.1]{panov-petrunin}, which we have been trying to prove for several years.

\begin{thm}{Conjecture}\label{braid-conj}
$\CB^n$ is $\CAT(0)$ for any integer $n\ge 1$.
\end{thm}

The case $n=1$ is trivial and the case $n=2$ is easy; $\CB^2$ is isometric to a product of the Euclidean plane and a Euclidean cone over the real line.
The case $n=3$ is already nontrivial;
two proofs are given by the authors \cite{panov-petrunin};
one of them follows from the result of Ruth Charney and Michael Davis \cite{charney-davis-93}.
Ruth Charney also proved a special case of the conjecture for $n=4$ \cite[Theorem~5.4]{charney2004}.
Namely, she proved that the real braid space with 4 strings is $\CAT(0)$.
One may think of the real braid space as the space of braid diagrams.
Formally, it is the inverse image of $\RR^n$ under the natural projection $\CB^n\to \CC^n$.
The real braid space is a convex subset of $\CB^n$. Hence it is $\CAT(0)$ if $\CB^n$ is $\CAT(0)$.

In this note we construct a class of convex functions on the braid space, which will be called \emph{pull-tight}.
The idea of this construction is inspired by Ian Agol's answer to our question \cite{agol}.

Let us describe pull-tight functions informally; an actual definition is given in Section~\ref{sec:pull-tight}.
Recall that $\WW^n$ is a configuration space of punctures in the plane.
Choose an isotopy class of a simple closed curve in the punctured plane and define a function $\tilde \WW^n\to\RR$ that returns the infimum of the lengths of curves in this class;
this function extends naturally to $\CB^n$.
There is also a variant of this function for simple open curves that connect two fixed ideal points of $\CC$.

\begin{thm}{Pull-tight theorem}\label{thm:pull-tight}
Every pull-tight function $\CB^n\to\RR$ is lower semicontinuous and convex.
\end{thm}

In the proof we need to show that $f\circ\gamma$ is convex for a pull-tight function $f$ and a geodesic $\gamma$ in $\CB^n$.
This can be easily done away from collisions of punctures along $\gamma$.
Assuming that convexity fails at a collision, we replace a segment of $\gamma$ near the collision by a shorter path with the same endpoints, and thus arrive at a contradiction.

\medskip

If $\CB^n$ is $\CAT(0)$, then it has no geodesic digons.
(Since $\CB^n$ is a polyhedral space and it has finitely many shapes, the converse holds as well; see \cite[II.5.4]{bridson-haefliger} and also compare with \cite[6.6]{AKP-invitation}.)

Recall that a \emph{geodesic loop} is a nonconstant closed curve that is locally geodesic away from its base point.
The pull-tight theorem implies the following statement, which is weaker than the absence of geodesic digons.

\begin{thm}{No-loops theorem}\label{thm:no-loops}
The braid space has no geodesic loops.
\end{thm}

In the proof we assume that a geodesic loop exists, and construct a pull-tight function that is nonnegative, not constant on the loop, and vanishes at the base point; this contradicts the pull-tight theorem.

\section{No-loops theorem}

\parbf{Proof of \ref{thm:no-loops} modulo \ref{thm:pull-tight}.}
Assume $\gamma\:[0,L]\to\CB^n$ is a geodesic loop;
denote by $p_1(t),\z\ldots,p_n(t)\in\CC$ its punctures at time $t$.
We may assume that
\[\mathrm{Re}(p_1(0))\le\ldots\le \mathrm{Re}(p_n(0)).\]
%Rotate the plane slightly so that points $p_i(0)$ have different real parts --- NO ---the puncture may coinside!!!

Choose $n-1$ pull-tight functions $f_1,\ldots,f_{n-1}$
such that the corresponding pull-tight curves are vertical straight lines that separate the punctures at time~$0$.
They may touch the punctures, but we remember on which side of each puncture the pull-tight curve passes;
the punctures can coincide and to make all this precise we have to use decorations defined in the next section.

Note that $f_i\ge0$ for each $i$.
Further, $f_i\circ\gamma(0)=f_i\circ\gamma(L)=0$;
since each $f_i$ is convex, $f_i\circ\gamma(t)\equiv 0$ for all~$i$.
It follows that the pull-tight lines form moving vertical barriers between the punctures, giving them no chance to interact; even if punctures meet or touch the barriers, they cannot go around each other.
In particular,
\[\mathrm{Re}(p_1(t))\le\ldots\le \mathrm{Re}(p_n(t))\]
for any $t$, and each puncture moves with constant velocity.
Thus the image of $\gamma$ in $\CC^n$ is a line segment --- a contradiction.
\qeds

\section{Pull-tight theorem}\label{sec:pull-tight}

\parbf{Collisions.}
Consider a geodesic $\gamma$ in $\CB^n$.
A point $\gamma(t)$ defines an array of punctures $p_1(t),\z\ldots,p_n(t)\in\CC$.
We say that $\gamma$ does not have a \emph{collision} at time $t$ if the velocities of all $p_i$ are constant in a neighborhood of $t$,
so the punctures move with constant velocity between the collision times.
(A priori, the set of collision times is merely closed.)

Since $\CB^n$ is a polyhedral space with finitely many shapes, the following lemma follows from the theorem of Martin R. Bridson and André Haefliger \cite[I.7.28]{bridson-haefliger}.
Formally, this theorem gives a bound on the number of collisions for geodesics of bounded length, but $\CB^n$ has a cone structure; therefore, we get a bound for arbitrary length.

\begin{thm}{Lemma}\label{lem:BFK}
There is a universal bound $N=N(n)$ on the number of collisions of geodesics in $\CB^n$.
\end{thm}

Let us mention that this lemma is also a limiting case of the result of Dmitri Burago, Serge Ferleger and Alexey Kononenko on the number of collisions in semi-dispersing billiards \cite{burago-ferleger-kononenko-1997,AKP-invitation};
their proof adapts to our situation and gets simpler.

\parbf{Construction.}
We will define a \emph{pull-tight} function $f\colon\CB^n\to\RR$ in two steps.
First, for a simple curve $\lambda$ in a punctured plane $\CC\setminus \{p_1,\ldots,p_n\}$ we define its \emph{pull-tight length}.
Then, using isotopies of $\lambda$ we get the pull-tight function on the whole space $\tilde \WW^n$ and extend it to $\CB^n$.

\parit{Pull-tight length.}
Choose an array of distinct punctures $\bm{p}=(p_1,\ldots,p_n)$ in $\CC$ and a simple oriented curve $\lambda$ in the punctured plane $\CC\setminus \{p_1,\ldots,p_n\}$. We call $\lambda$ the \emph{defining curve}.
The curve can be closed or open;
in the latter case, we assume that $\lambda$ connects two fixed ideal points of $\CC$.

{

\begin{wrapfigure}{o}{40mm}
\centering
\vskip-2mm
\includegraphics{mppics/pic-70}
\vskip-0mm
\end{wrapfigure}

Think of $\lambda$ as a thread and make it tight in the punctured plane.
The tightened curve is a polygonal line, say $\ell$, with vertices at the punctures;
self-intersections are possible.
If $\lambda$ is open, then $\ell$ has two infinite edges, or might be a straight line without vertices.
The polygonal line $\ell$ can be arbitrarily well approximated by a simple curve isotopic to~$\lambda$.
We will refer to $\ell$ as the \emph{pull-tight curve}.

}

If $\lambda$ is closed, define the \emph{pull-tight length} as the length of $\ell$; this is the infimum length over the isotopy class of $\lambda$.

If $\lambda$ is open,  $\ell$ might not exist.
This happens if $\lambda$ can be pushed into a strictly convex subset of $\CC$ containing no punctures, so it is possible to shrink the thread indefinitely.
In this case, we say that $\lambda$ is an \emph{escaping} curve and assume that the pull-tight length is undefined (one could also define it as $-\infty$).
Otherwise the pull-tight length is defined as the excess of its length compared to the curve composed of two rays emanating from the origin.
More formally, suppose that $\alpha$ and $\omega$ are the ideal points at the beginning and the end,
and $p_{i_1},\ldots,p_{i_k}$ is the sequence of vertices of $\ell$.
Then the pull-tight length is
\[|\alpha-p_{i_1}|^\flat+|p_{i_1}-p_{i_2}|+\ldots+|p_{i_{k-1}}-p_{i_k}|+|p_{i_k}-\omega|^\flat,\]
where $|\xi-x|^\flat$
denotes the Busemann function at $x$ with ideal point $\xi$ normalized so that $|\xi-0|^\flat=0$; thus,
\[|\xi-x|^\flat=-\langle u_\xi,x\rangle,\]
where $u_\xi$ denotes the unit vector in the direction of $\xi$.
Note that $|\xi-x|^\flat$ might take negative values. If $\ell$ has no vertices, then it is a straight line and the pull-tight length is defined to be zero.

\parit{Pull-tight function on $\tilde \WW^n$ and $\CB^n$.}
Let again $\bm p=(p_1,\ldots,p_n)$ be a fixed array of distinct punctures in $\CC$ and $\lambda$ an oriented simple nonescaping curve in the punctured plane $\CC \setminus \{p_1,\ldots,p_n\}$.
A point $\bm{\tilde q}$ of $\tilde \WW^n$ corresponds to another array of punctures $\bm q=(q_1,\ldots,q_n)\in \WW^n$ and a homotopy class of paths in $\WW^n$ connecting $\bm p$ and $\bm q$.
Thus we can isotope $\lambda$ to a simple curve $\lambda'$ lying in the punctured plane $\CC\setminus \{q_1,\ldots,q_n\}$.
Define $f(\bm{\tilde q})$ as the pull-tight length of $\lambda'$.
Thus we get a function $f$ on $\tilde \WW^n$.

Note that $f$ is locally Lipschitz, and it naturally extends to a lower semicontinuous function on $\CB^n$.
Along a geodesic this extension is continuous: at a collision the edges inside the colliding clusters shrink to zero, while all remaining terms have their usual limits.

\parbf{Proof of \ref{thm:pull-tight}.}
Let $f\:\CB^n\to \RR$ be a pull-tight function with defining curve $\lambda$.
Given a unit-speed geodesic $\gamma$ in $\CB^n$,
we need to show that $f\circ\gamma$ is convex.

Each point $\gamma(t)$ defines an array of punctures $p_1(t),\ldots,p_n(t)\in\CC$.
The curve $\lambda$ divides the punctures into two sets that lie on its left and right sides.

\parit{Encoding.}
Denote by $\ell(t)$ the pull-tight curve at time $t$.
By construction, $\ell(t)$ is a polygonal line with vertices only at $p_i(t)$,
and $\ell(t)$ can be encoded by its vertex sequence --- the ordered list of its vertex labels, with repetitions allowed.
(For an open curve, an empty sequence corresponds to a straight line.)
If the curve passes a puncture without turning, we do not record that occurrence as a vertex.

\parit{Decoration.}
By construction, $\ell(t)$ and $p_i(t)$ can be perturbed so that punctures become distinct and the curve becomes simple and isotopic to $\lambda$;
this can be done continuously and simultaneously for all $t$.
The isotopy information carried by these perturbations will be called the \emph{decoration} of $\gamma$ and the pull-tight curve.

\parit{Changing encoding.}
Suppose that the vertex sequence of $\ell(t)$ does not change in some open time interval $\II$ and no collision occurs in $\II$.
Note that if two points move with constant velocities in the plane, then the distance between them is a convex function of time.
It follows that $f\circ\gamma|_{\II}$ is convex.

The vertex sequence can change in several ways:
(1) when a puncture hits or leaves $\ell(t)$,
(2)  at a collision or an isolated coincidence of punctures,
and (3) in the mixed case, when such events happen simultaneously.
Such events happen at finitely many times, say $t_1,\ldots,t_m$.
Indeed, by \ref{lem:BFK}, we have only finitely many collisions.
Between the collision times the velocities of punctures are constant;
any two punctures either meet at most once or coincide thruout the interval, and hitting/leaving events are described by isolated zeros of algebraic equations.

For each $t_i$, we have to check that
\[(f\circ\gamma)'(t_i+) - (f\circ\gamma)'(t_i-)\ge 0;
\eqlbl{1star}\]
here $(f\circ\gamma)'(t\pm)$ denote one-sided derivatives.
Once this is done, convexity of $f\circ\gamma$ follows.

\parit{Hits and leaves.}
Suppose that the vertex sequence changes at time $t_i$  without a collision or an isolated coincidence of punctures.
For example, $\ell(t)$ catches up with the point $p_5(t)$, so $p_5(t)$ enters the vertex sequence between vertices $p_2(t)$ and $p_7(t)$, as in the picture.

\begin{wrapfigure}{o}{46mm}
\centering
\vskip-0mm
\includegraphics{mppics/pic-50}
\vskip-0mm
\end{wrapfigure}

In a neighborhood of $t_i$ the curve is described by two vertex sequences;
denote the corresponding functions by $f_{\textrm{before}}$ and $f_{\textrm{after}}$;
if $t\approx t_i$, then
\[f\circ\gamma(t)=
\begin{cases}
f_{\textrm{before}}\circ\gamma(t)&\text{if\ } t\le t_i,
\\
f_{\textrm{after}}\circ\gamma(t)&\text{if\ } t\ge t_i.
\end{cases}
\]
Note that
\[(f_{\textrm{before}}\circ\gamma)'(t_i)=(f_{\textrm{after}}\circ\gamma)'(t_i);\]
hence \ref{1star} follows.

\parit{Rescaling.}
Suppose that $k$ punctures coincide at time $t_i$; let us keep calling it a \emph{collision}.
We can assume that $t_i=0$ and the collision (or coincidence) point is also $0\in\CC$.

Rescale a neighborhood of the collision point and a time interval around the collision time, keeping the velocities the same.
In the limit we get a collision of all $k$ punctures in $\CB^k$.
The punctures come with constant velocities to $0$ and run from it with constant velocities.
Denote by $\tilde \ell(t)$, $\tilde \gamma$, and $\tilde p_1,\ldots,\tilde p_k$
the limits of the pull-tight curve, the geodesic in $\CB^k$, and the corresponding punctures.

Suppose that $\tilde \ell(t)$ is still defined by the isotopy class of one simple curve $\tilde\lambda$ in the $k$-punctured plane with fixed ideal endpoints.
Let us check \ref{1star} for the corresponding pull-tight function $\tilde f\colon\CB^k\to \RR$.

Note that the shape of $\tilde \ell(t)$ changes only at $t=0$:
for $t>0$ and $t<0$ it only rescales.
At $t=0$ these shapes approach a plane angle from both sides;
let $\theta$ be its angle measure.
If \ref{1star} does not hold, then $\theta<\pi$.
In other words, the ideal points of $\tilde \ell$ cannot be opposite; if they are, then $(\tilde f\circ\tilde \gamma)(t)\ge 0$ for any $t$ and $(\tilde f\circ\tilde \gamma)(0)= 0$, and thus \ref{1star} holds for $\tilde f\circ\tilde \gamma$.

Since $\theta<\pi$, we can add a constant vector to all velocities in the $k$-puncture array and reduce the problem to the case
\[(\tilde f\circ\tilde \gamma)'(0-)>0>(\tilde f\circ\tilde \gamma)'(0+).\]
We may still assume that the collision occurs at the origin at time zero;
so, $\tilde f\circ\tilde \gamma$ is a piecewise linear function, with exactly one kink at time zero,
$\tilde f\circ\tilde \gamma(0)\z=0$, and, for some $\eps>0$,
\[\tilde f\circ\tilde \gamma(t)\le -\eps\cdot |t|.
\eqlbl{sharp-max}\]
We are going to construct a shortcut of $\tilde \gamma$, showing that it is not geodesic.

\begin{figure}[b!]
\centering
\vskip-0mm
\includegraphics{mppics/pic-60}
\vskip-0mm
\end{figure}

\parit{$\CAT(0)$ sides.}
The curve $\tilde\ell(t)$ might have self-intersections, but it admits an approximation by a simple curve, which divides the plane into two regions on its left and right sides.
Let us equip these regions with the induced length metric and pass to their completions;
we get $\CAT(0)$ spaces \cite[7.12]{AKP-invitation}.
Pass to their pointed limits as the approximating simple curve approaches $\tilde\ell(t)$; denote the limit spaces by $L_t$ and $R_t$; they are $\CAT(0)$ as limits of $\CAT(0)$ spaces and $\tilde\ell(t)$ runs in their metric boundaries.
(If $\tilde\ell(t)$ is simple, then $L_t$ and $R_t$ are just the closed regions on the left and right sides of $\tilde\ell(t)$; the described construction generalizes this picture.)

One of $L_0$ and $R_0$ is a convex infinite sector of angle $\theta$;
let us assume that this is~$L_0$.
Choose an ideal point $\xi$ that lies midway between the ideal points, say $\alpha$ and $\omega$, of $\tilde\ell(t)$.
Let us denote by $x\mapsto|x-\xi|^\flat_t$ its Busemann function $L_t\to\RR$.
These functions are defined in the length metric of $L_t$.
Formally speaking, $\xi$ corresponds to an ideal point in each space $L_t$,
but far from $0$, the half-lines in $L_t$ in the direction of $\xi$ become parallel.
Thus the Busemann functions $x\z\mapsto |x-\xi|^\flat_t$ can be normalized uniformly.
(Far away from $0$, they all agree with the Busemann function for $\xi$ with respect to the Euclidean metric: $x\mapsto-\langle x,u_\xi\rangle$.)
In particular we have $|t\cdot x-\xi|^\flat_t=t\cdot |x-\xi|^\flat_1$ for any $t\ge 0$ and $x\in L_1$.
The superlevel sets of this Busemann function are compact sets in $L_t$.
Moreover,
\[2\cdot \cos \tfrac\theta2 \cdot |x-\xi|^\flat_t
\le
|\alpha-x|^\flat_t+|x-\omega|^\flat_t
\le \tilde f\circ \tilde\gamma(t)\]
for any $x\in L_t$.
The first inequality follows from $\CAT(0)$ comparison and the second from the triangle inequality, both applied to $x$ and faraway points in the directions of $\xi$, $\alpha$, and $\omega$.
More precisely, the first inequality is the inequality between the median and an equal side of a huge isosceles triangle (in the Euclidean plane, it becomes an identity).
To prove the second inequality, note that $|\alpha-x|^\flat_t+|x-\omega|^\flat_t$ can be interpreted as the right-hand side after we replace all punctures on the left of the curve by one puncture at $x$.
This changes the isotopy class of admissible curves, but since the old minimizer remains admissible, this can only reduce the minimal length.

Since $\tilde f\circ\tilde \gamma(t)\le -\eps\cdot |t|$, we get
\[|x-\xi|^\flat_t\le -\eps'\cdot|t|\]
for any $x\in L_t$ and some fixed $\eps'>0$.

Suppose $\tilde p_i$ lies on the left side of~$\tilde \lambda$.
Let $\delta_i=-|\tilde p_i(1)-\xi|^\flat_1$;
note that $\tilde p_i(t)\z=t\cdot \tilde p_i(1)$ and $|\tilde p_i(t)-\xi|^\flat_t=-t\cdot \delta_i$ for any $t\ge 0$.
Let us connect $\tilde p_i(1)$ to $\xi$ by a unit-speed geodesic half-line $\sigma_i\:[\delta_i,\infty)\z\to L_1$;
so $|\sigma_i(s)-\xi|^\flat_1+s\equiv0$.
For any $t>0$,
\[s\mapsto t\cdot \sigma_i(s/t)
\eqlbl{2star}\]
defines the unit-speed geodesic half-line in $L_t$ from $\tilde p_i(t)$ to $\xi$.
Further, if $s$ is large, then
\[\sigma_i(s)=s\cdot u_\xi+ z_i,\eqlbl{3star}\]
for a complex constant $z_i$ such that $z_i\perp u_\xi$.

\parit{Pushing the punctures.}
We now push the punctures on the left toward the bisector direction $\xi$, stopping them at a fixed Busemann level, and show that the resulting path has the same endpoints but smaller energy and therefore length.

Choose $\delta>0$ such that $\delta<\delta_i$ for any $\tilde p_i$ on the left side of $\tilde \lambda$.
Let $\hat p_i(t)$ be $\tilde p_i(t)$ moved toward $\xi$ along the half-line defined by \ref{2star}:
if $\tilde p_i(t)$ is below the $(-\delta)$-level of the Busemann function for $\xi$, then $\tilde p_i(t)$ is not moved;
otherwise it is moved to the $(-\delta)$-level.
In other words, \[\hat p_i(t)=t\cdot \sigma_i(\max\{\,\delta/t, \delta_i\,\});\]
so, if $t\cdot \delta_i\le \delta$, then
we move $\tilde p_i(1)$ along $\sigma_i$ to the $(-\delta/t)$-level and then scale it down by $t$.

Note that $|\tilde p_i'(t)|^2=|\tilde p_i(1)|^2$ for any $t>0$.
Therefore
\[\int\limits_{0}^{\delta/\delta_i}|\tilde p_i'(t)|^2\cdot dt
=
(\delta/\delta_i)\cdot|\tilde p_i(1)|^2.\]
If $0<t<\delta/\delta_i$, then $\hat p_i(t)=t\cdot \sigma_i(\delta/t)$
and
\[|\sigma'_i|=1
\quad\Rightarrow\quad
|\hat p_i'(t)|^2=\left(\frac{|\hat p_i(t)|^2-\delta^2}{t}\right)'\]
almost everywhere.
By \ref{3star}, we have $|\hat p_i(t)|^2-\delta^2=o(t)$, whence
\[\int\limits_0^{\delta/\delta_i}|\hat p_i'(t)|^2\cdot dt
=
(\delta/\delta_i)\cdot(|\tilde p_i(1)|^2-\delta_i^2)
<
\int\limits_0^{\delta/\delta_i}|\tilde p_i'(t)|^2\cdot dt.\]

Observe that all modified punctures coincide at time $0$: if $\tilde p_i$ lies on the left side of $\tilde\lambda$, then $\hat p_i(0)=\delta\cdot u_\xi$.
Therefore, repeating the same construction for $t\le 0$ produces a curve $t\mapsto \hat p_i(t)$ defined on the whole $\RR$.

Each puncture $\tilde p_i$ on the right side will not be modified;
that is, we let $\hat p_i(t)=\tilde p_i(t)$.
Thus we obtain a new curve $t\mapsto \hat p_i(t)$ for each $i$.
These curves correspond to a new curve $\hat\gamma$ in $\CB^k$ that coincides with $\tilde \gamma$ away from some time interval $[-R,R]$.
Indeed, punctures cannot cross $\tilde\ell(t)$; so the punctures on its different sides do not interact.
On the right side the punctures move as before,
and all modified punctures approach $\delta\cdot u_\xi$ at time $0$,
so the decoration within this colliding cluster can be chosen to match the original braid outside $[-R,R]$.

The above inequality implies that
\[\int\limits_{-R}^{R}|\tilde p'_i(t)|^2\cdot dt
\ge
\int\limits_{-R}^{R}|\hat p'_i(t)|^2\cdot dt\]
for each $i$; moreover, if $\tilde p_i$ is on the left side then the inequality is strict.
Since $\tilde\gamma$ has constant speed,
summing over $i$ and applying Cauchy--Bunyakovsky, we obtain
\[\length \tilde\gamma|_{[-R,R]}>\length\hat\gamma|_{[-R,R]},\]
which gives the required contradiction.

The limit curve $\tilde\ell(t)$ might be empty, but then the collision makes no contribution to the derivative jump of $f\circ\gamma$.
Further, $\tilde\ell(t)$ might be closed, but then $\tilde f\circ \tilde\gamma$ cannot satisfy \ref{sharp-max} (in this case $\tilde f\circ \tilde\gamma(t)\ge0$ and $\tilde f\circ \tilde\gamma(0)=0$).
Furthermore, $\tilde\ell(t)$ might have several components;
they have to be described by isotopy classes of nonintersecting simple curves $\tilde\lambda_1,\ldots,\tilde\lambda_m$ with fixed ideal endpoints.
Denote by $\tilde f_1,\dots,\tilde f_m$ the corresponding pull-tight functions.
One of the functions $\tilde f_i\circ\tilde\gamma$ must have a negative derivative jump at $0$ (the original derivative jump is the sum of the component jumps).
Applying the common-velocity adjustment, we can assume that $\tilde f_i\circ\tilde\gamma$
has a sharp maximum at zero; in other words, $\tilde f_i\circ\tilde\gamma$ satisfies \ref{sharp-max}.
Then we can repeat the argument for $\tilde\lambda_i$.
As before, we get a shortcut $\hat\gamma$ for $\tilde \gamma$ and arrive at a contradiction.

\parit{Mixed case.} Note that our calculations are local: combining them, one can prove \ref{1star} in the case when several events (collisions, isolated coincidences, and hitting/leaving) happen simultaneously.
\qeds

\section{Postscript}

The following observations might save time for someone interested in the conjecture;
this section explains where we hit the wall.
We still have a feeling that one new geometric idea can finish the story.

In the definition of $\CB^n$ one can lift an arbitrary diagonal metric on $\CC^n$.
In other words, we prescribe positive masses $\mu_1,\z\ldots,\mu_n$ to the punctures $p_1(t),\z\ldots, p_n(t)$ that move in $\CC$ and define the speed of the corresponding path in $\CB^n$ as
\[\sqrt{\mu_1\cdot|p_1'(t)|^2+\ldots+ \mu_n\cdot|p_n'(t)|^2}.\]
Let us denote the obtained space by $\CB^n(\bm{\mu})$, where $\bm{\mu}=(\mu_1,\ldots,\mu_n)$.%
\footnote{Equivalently, we can apply a diagonal matrix with $\sqrt{\mu_1},\ldots,\sqrt{\mu_n}$ on the diagonal to the complex hyperplane arrangement $z_i=z_j$ (so the new complex arrangement is defined by $\frac{z_i}{\sqrt{\mu_i}}=\frac{z_j}{\sqrt{\mu_j}}$).
Then $\CB^n(\bm{\mu})$ is isometric to the universal cover of $\CC^n$ (with the standard metric) infinitely branching along the new hyperplane arrangement.}

The definition of a pull-tight function extends to $\CB^n(\bm{\mu})$.
Moreover, the proof of the pull-tight theorem has a straightforward generalization that gives the following.

\begin{thm}{Observation}
Every pull-tight function $\CB^n(\bm{\mu})\to\RR$ is lower semicontinuous and convex.
In particular, $\CB^n(\bm{\mu})$ does not have geodesic loops.
\end{thm}

\begin{thm}{Corollary from the conjecture}\label{cor:masses}
If the conjecture holds, then $\CB^n(\bm{\mu})$ is $\CAT(0)$ for any array of positive masses $\bm{\mu}=(\mu_1,\ldots,\mu_n)$.
\end{thm}

\parbf{Proof.}
If the masses $\mu_i$ are integers, then instead of the $i$-th string of the braid one can consider $\mu_i$ identical copies of this string.
Thus we can think that $\CB^n(\bm{\mu})$ is a subset of $\CB^M$ with $M=\mu_1+\ldots+\mu_n$.

Note that $\CB^n(\bm{\mu})$ is convex in $\CB^M$.
Indeed, forgetting all but one copy of each of $n$ groups of strings defines $\mu_1\cdots\mu_n$ projections $\pi\:\CB^M\to \CB^n(\bm{\mu})$; every projection restricts to the identity on $\CB^n(\bm{\mu})\subset\CB^M$.
If $\gamma$ has endpoints in $\CB^n(\bm{\mu})$, then every projected path has the same endpoints.
The energy $E(\gamma)$ equals the average of the energies $E(\pi\circ\gamma)$ for all projections $\pi$.
Therefore  $E(\gamma)\ge E(\pi\circ\gamma)$ for some projection~$\pi$.
We can assume that $\gamma$ is parametrized with constant speed; in this case, \[\length\gamma\ge \length(\pi\circ\gamma).\]
Applying the inequality to a minimizing ambient geodesic $\gamma$ establishes weak convexity.
Uniqueness under the conjecture then gives convexity.

In particular, the conjecture for $\CB^M$ implies the corollary for $\CB^n(\bm{\mu})$.
By rescaling, rational approximation, and passing to the limit we get the case of arbitrary positive masses.
\qeds

\parbf{Approach to the conjecture.} Let us try to prove \ref{cor:masses} by induction on $n$.
The case $n=3$ already follows from \cite{charney-davis-93} (checking the angle condition is straightforward).
The first new case $n=4$ might go as follows.

Let us split off the diagonal translation factor $\CC$; in other words, we assume that the center of mass of the punctures is at $0$.
The remaining factor is a metric cone; denote its unit link by $\Sigma^4(\bm{\mu})$.
We need to show that it is $\CAT(1)$.

By the induction hypothesis $\Sigma^4(\bm{\mu})$ is locally $\CAT(1)$.%
\footnote{\label{foot:diagonal} In fact our theorem implies that every $\tfrac\pi2$-ball in $\Sigma^4(\bm{\mu})$ equipped with the induced length metric is $\CAT(1)$;
in particular, $\Sigma^4(\bm{\mu})$ is $\CAT(4)$ and any $\tfrac\pi4$-ball in $\Sigma^4(\bm{\mu})$ is $\CAT(1)$.}
Assume $\Sigma^4(\bm{\mu})$ is not $\CAT(1)$.
By \cite[II.4.17]{bridson-haefliger}, it contains a shortest geodesic circle, say $\zeta$, of length $L(\bm{\mu})<2\cdot\pi$ (a \emph{geodesic circle} is a distance-preserving embedding of a circle with its length metric).
If one of the masses approaches zero, we expect to have $L(\bm{\mu})\z\to2\cdot \pi$ by the induction hypothesis.

Any pull-tight function $f$ is spherically convex on $\zeta$;
that is $(f\circ \zeta)''+f\circ \zeta\ge 0$ in the barrier sense.
Thus if we can find a pull-tight function $f$ (or a lasso function, see below) that is nonpositive on two opposite points of $\zeta$ and $f\circ\zeta\ne0$, then the problem will be solved.
This approach might work in principle, but it requires geometric information on $\zeta$:
\textsc{there are closed curves in $\Sigma^4(\bm{\mu})$ of length $<2\cdot \pi$ that do not meet this requirement for any choice of pull-tight function}.

Given an array of punctures $p_1,\ldots, p_4$, we translate it to make its center of mass zero, then normalize it so that $\mu_1\cdot|p_1|^2+\ldots+ \mu_4\cdot|p_4|^2=1$;
the resulting configuration corresponds to a point in $\Sigma^4(\bm{\mu})$.
This makes it possible to identify the underlying sets of $\Sigma^4(\bm{\mu})$ for different masses; so we may think of $\Sigma^4(\bm{\mu})$ as the same set equipped with different metrics.

The geodesic circle $\zeta$ is not null-homotopic in the class of curves of length $<2\cdot\pi$;
that is, it cannot be contracted thru curves all of length $<2\cdot\pi$.
After a small perturbation of the masses, $\zeta$ may no longer be a geodesic,
but shrinking its length must produce a closed geodesic $\zeta'$, which cannot be trivial.
This gives an upper bound on $L(\bm{\mu})$ after small perturbation of masses.

Since the rescaling of masses does not change the isometry type of $\Sigma^4(\bm{\mu})$, we get
\[\mu_1\cdot \frac{\partial L}{\partial \mu_1}+\ldots+\mu_4\cdot \frac{\partial L}{\partial \mu_4}=0\] if the left-hand side is defined.

The function $\bm{\mu}\mapsto L(\bm{\mu})$ need not be smooth,
but it is semiconcave.
At a minimum point of $L$ it must be differentiable (the induction hypothesis should imply that the minimum is attained in the interior), and we have $\frac{\partial L}{\partial \mu_i}=0$ for every $i$.
Assume that a unit-speed parametrization of the geodesic circle is given by punctures $p_1(t),\ldots,p_4(t)$; so
\begin{align*}
\mu_1\cdot p_1(t)+&\ldots+ \mu_4\cdot p_4(t)=0,
\\ \mu_1\cdot|p_1(t)|^2+&\ldots+ \mu_4\cdot|p_4(t)|^2=1,
\\ \mu_1\cdot|p_1'(t)|^2+&\ldots+ \mu_4\cdot|p_4'(t)|^2=1.
\end{align*}
Moreover, $p_i''(t)+p_i(t)=0$ if $p_i$ does not collide at time $t$.
The condition $\frac{\partial L}{\partial \mu_i}=0$ means that
\[\sum_{t_j}\langle\, p_i(t_j),p_i'(t_j+)-p_i'(t_j-)\, \rangle=0,\]
where the sum is taken at the collision times~$t_j$.

The last identity might contradict $L<2\cdot\pi$.
Also, the estimate on $L$, which is described above, might imply that the Hessian of $L$ is negative (in the barrier sense) in some direction, which should prove the conjecture.
\textsc{The calculations seem to be simpler for the Fisher metric on the mass simplex:} $g(h,h)=\tfrac{h_1^2}{\mu_1}+\ldots+\tfrac{h_4^2}{\mu_4}$, where $h_1+\ldots+h_4=0$; compare to the description of $\CB^4(\bm{\mu})$ in the footnote \ref{foot:diagonal}.

Note, however, that for an equator the Hessian and Laplacian vanish.
Therefore, one has to use that $L<2\cdot\pi$ and the geometry of the geodesic circle.
The following observation might help: since $L<2\cdot\pi$, the image of the geodesic circle in $\mathbb{S}^5\subset \CC^3$ lies in a hemisphere.

\parbf{Lasso functions.}
The proof of convexity of pull-tight functions uses very little and admits variations;
here is one example.

Choose one puncture $p_i$, called the \emph{caught puncture}.
Pass to the covering of $\CC$ infinitely branching at all punctures except $p_i$, and lift $p_i$ to this covering.
Consider the Busemann function associated with the chosen ideal point $\tilde \xi$ of the covering, evaluated at this lift.
Our proof of convexity works for this function as well.
The obtained function $s$ (let's call it a \emph{lasso function})
has a closely related pull-tight function $f$ --- the one with both ideal endpoints at the projection $\xi$ of $\tilde \xi$ and $\lambda$ surrounding just $p_i$.

For the compatible normalizations we have $2\cdot s\le f$; in particular, the subzero set of $s$ contains that of $f$.
So, in a sense $s$ looks more like a Busemann function on $\CB^n$.
(If Busemann functions could be approximated by convex functions, then the conjecture would follow, but one may need much less.)

\parbf{Acknowledgements.}
The postscript contains two essential suggestions from AI (GPT-6 Astra); both are marked in \textsc{small caps}.

{\sloppy
\def\emph{\textit}
\printbibliography[heading=bibintoc]

@book{AKP-invitation,
    AUTHOR = {Alexander, S. and Kapovitch, V. and Petrunin, A.},
     TITLE = {An invitation to Alexandrov geometry: CAT(0) spaces},
    SERIES = {Lecture Notes in Mathematics},
    VOLUME = {2395},
   EDITION = {2},
 %PUBLISHER = {Springer, Cham},
      YEAR = {2026},
    % PAGES = {xv+124},
      ISBN = {978-3-032-25348-4; 978-3-032-25349-1},
   MRCLASS = {99-06},
  MRNUMBER = {5127659},
       DOI = {10.1007/978-3-032-25349-1},
       URL = {https://doi.org/10.1007/978-3-032-25349-1},
}

@article {allcock,
    AUTHOR = {Allcock, D.},
     TITLE = {Completions, branched covers, {A}rtin groups, and singularity
              theory},
   JOURNAL = {Duke Math. J.},
  FJOURNAL = {Duke Mathematical Journal},
    VOLUME = {162},
      YEAR = {2013},
    NUMBER = {14},
     PAGES = {2645--2689},
      ISSN = {0012-7094,1547-7398},
   MRCLASS = {53C23 (14B07 20F36 51K10 57N65)},
  MRNUMBER = {3127810},
MRREVIEWER = {Igor\ Belegradek},
       DOI = {10.1215/00127094-2380977},
       URL = {https://doi.org/10.1215/00127094-2380977},
}

@book {bridson-haefliger,
    AUTHOR = {Bridson, M. and Haefliger, A.},
     TITLE = {Metric spaces of non-positive curvature},
    SERIES = {Grundlehren der mathematischen Wissenschaften},
    VOLUME = {319},
% PUBLISHER = {Springer-Verlag, Berlin},
      YEAR = {1999},
     %PAGES = {xxii+643},
      ISBN = {3-540-64324-9},
   MRCLASS = {53C23 (20F65 53C70 57M07)},
  MRNUMBER = {1744486},
MRREVIEWER = {Athanase\ Papadopoulos},
       DOI = {10.1007/978-3-662-12494-9},
       URL = {https://doi.org/10.1007/978-3-662-12494-9},
}

@article {burago-ferleger-kononenko-1997,
    AUTHOR = {Burago, D. and Ferleger, S. and Kononenko, A.},
     TITLE = {Uniform estimates on the number of collisions in
              semi-dispersing billiards},
   JOURNAL = {Ann. of Math. (2)},
  FJOURNAL = {Annals of Mathematics. Second Series},
    VOLUME = {147},
      YEAR = {1998},
    NUMBER = {3},
     PAGES = {695--708},
      ISSN = {0003-486X},
   MRCLASS = {58F11 (58F99 70F10 70F35)},
  MRNUMBER = {1637663},
MRREVIEWER = {Nikolai Chernov},
       DOI = {10.2307/120962},
       URL = {https://doi.org/10.2307/120962},
}

@article {charney-davis-93,
    AUTHOR = {Charney, R. and Davis, M.},
     TITLE = {Singular metrics of nonpositive curvature on branched covers
              of {R}iemannian manifolds},
   JOURNAL = {Amer. J. Math.},
  FJOURNAL = {American Journal of Mathematics},
    VOLUME = {115},
      YEAR = {1993},
    NUMBER = {5},
     PAGES = {929--1009},
      ISSN = {0002-9327,1080-6377},
   MRCLASS = {53C23 (57M12)},
  MRNUMBER = {1246182},
MRREVIEWER = {Athanase\ Papadopoulos},
       DOI = {10.2307/2375063},
       URL = {https://doi.org/10.2307/2375063},
}

@article {charney-davis-95,
    AUTHOR = {Charney, R. and Davis, M.},
     TITLE = {The {$K(\pi,1)$}-problem for hyperplane complements associated
              to infinite reflection groups},
   JOURNAL = {J. Amer. Math. Soc.},
  FJOURNAL = {Journal of the American Mathematical Society},
    VOLUME = {8},
      YEAR = {1995},
    NUMBER = {3},
     PAGES = {597--627},
      ISSN = {0894-0347,1088-6834},
   MRCLASS = {52B30 (20F55 55P20)},
  MRNUMBER = {1303028},
MRREVIEWER = {Hiroaki\ Terao},
       DOI = {10.2307/2152924},
       URL = {https://doi.org/10.2307/2152924},
}

@article {charney2004,
    AUTHOR = {Charney, R.},
     TITLE = {The {D}eligne complex for the four-strand braid group},
   JOURNAL = {Trans. Amer. Math. Soc.},
  FJOURNAL = {Transactions of the American Mathematical Society},
    VOLUME = {356},
      YEAR = {2004},
    NUMBER = {10},
     PAGES = {3881--3897},
      ISSN = {0002-9947,1088-6850},
   MRCLASS = {20F36 (20F55 32S22 52C35 57M07)},
  MRNUMBER = {2058510},
MRREVIEWER = {Stephen\ P.\ Humphries},
       DOI = {10.1090/S0002-9947-03-03425-1},
       URL = {https://doi.org/10.1090/S0002-9947-03-03425-1},
}

@article {panov-petrunin,
    AUTHOR = {Panov, D. and Petrunin, A.},
     TITLE = {Ramification conjecture and {H}irzebruch's property of line
              arrangements},
   JOURNAL = {Compos. Math.},
  FJOURNAL = {Compositio Mathematica},
    VOLUME = {152},
      YEAR = {2016},
    NUMBER = {12},
     PAGES = {2443--2460},
      ISSN = {0010-437X,1570-5846},
   MRCLASS = {53C23 (52C35 57N65)},
  MRNUMBER = {3594282},
MRREVIEWER = {Alexander\ O.\ Ivanov},
       DOI = {10.1112/S0010437X16007648},
       URL = {https://doi.org/10.1112/S0010437X16007648},
}

@article {paris,
    AUTHOR = {Paris, L.},
     TITLE = {{$K(\pi,1)$} conjecture for {A}rtin groups},
   JOURNAL = {Ann. Fac. Sci. Toulouse Math. (6)},
  FJOURNAL = {Annales de la Facult\'e{} des Sciences de Toulouse.
              Math\'ematiques. S\'erie 6},
    VOLUME = {23},
      YEAR = {2014},
    NUMBER = {2},
     PAGES = {361--415},
      ISSN = {0240-2963,2258-7519},
   MRCLASS = {57M07 (20F36 20F55 52C35 55P20)},
  MRNUMBER = {3205598},
MRREVIEWER = {Jon\ McCammond},
       DOI = {10.5802/afst.1411},
       URL = {https://doi.org/10.5802/afst.1411},
}

@misc{petrunin2026,
      title={A gallery of open problems in geometry that keep me puzzled and amused},
      author={Anton Petrunin},
      year={2026},
      eprint={2510.11908},
      archivePrefix={arXiv},
      primaryClass={math.DG},
      url={https://arxiv.org/abs/2510.11908},
}

@phdthesis{vanderlek,
  author = {van der Lek, H.},
  title = {The Homotopy Type of Complex Hyperplane Complements},
  school = {University of Nijmegen},
  year = {1983},
}

@MISC {agol,
    TITLE = {Braid group on 4 strands},
    AUTHOR = {I. Agol},
    HOWPUBLISHED = {MathOverflow},
    %NOTE = {URL:https://mathoverflow.net/q/307215 (version: 2018-07-31)},
    EPRINT = {https://mathoverflow.net/q/307215},
    URL = {https://mathoverflow.net/q/307215}
}
\fussy
}

\end{document}